\documentclass[11pt]{article}

\usepackage{soul}
\usepackage{lipsum}
\usepackage{amsfonts}
\usepackage{amscd}
\usepackage{amssymb}
\usepackage{amsthm}
\usepackage{amsmath, xspace}
\usepackage{blkarray}
\usepackage{cancel}
\usepackage{graphics}
\usepackage{graphicx}
\usepackage{enumitem}
\usepackage{fancyhdr}
\usepackage{mathdots}
\usepackage{mathrsfs}
\usepackage{multicol}
\usepackage{stmaryrd}
\usepackage{ytableau}

\usepackage[plainpages,backref]{hyperref}

\usepackage{lscape}

\theoremstyle{plain}
\newtheorem{thm}{Theorem}[section]

\theoremstyle{definition}

\newtheorem{example}{Example}[section]

\theoremstyle{remark}

\numberwithin{equation}{section}

\providecommand{\keywords}[1]{\textbf{\textit{Key words---}} #1}

\def\ZZ{\mathbb{Z}}

\def\id{\mathrm{id}}

\def\wt{\mathrm{wt}}

\usepackage{array}

\newsavebox\topalignbox
\newcolumntype{T}{%
  >{\begin{lrbox}\topalignbox
    $}
  c
  <{\!$\end{lrbox}%
    \raisebox{\dimexpr-\height+\ht\strutbox\relax}%
      {\usebox\topalignbox}}}
\newcolumntype{M}{>{$\vcenter\bgroup\hbox\bgroup}c<{\egroup\egroup$}}

\usepackage{color}
\definecolor{dgreen}{rgb}{0, .65, 0.15}
\definecolor{dred}{rgb}{.7, .2, 0}
\definecolor{dblue}{rgb}{.3,0,.7}
\definecolor{dpurple}{rgb}{.5,0,.65}

\newcommand{\blue}[1]{{\color{blue}#1}}
\newcommand{\red}[1]{{\color{dred}#1}}

\usepackage{tikz}
	\usepgflibrary{shapes.geometric}
	\usetikzlibrary{shapes.misc}

\tikzstyle{V}=[draw, fill =black, circle, inner sep=0pt, minimum size=1.5pt]
\tikzstyle{wV}=[draw, fill =white, circle, inner sep=0pt, minimum size=4.5pt]
\tikzstyle{bV}=[draw, fill =black, circle, inner sep=0pt, minimum size=4.5pt]

\newcommand{\TikZ}[1]{
\begin{matrix}\begin{tikzpicture}#1\end{tikzpicture}\end{matrix}
}

\newcounter{bx}
\newcounter{r}
\newcounter{s}
\newcommand\Part[1]{
        \setcounter{r}{1}
        \setcounter{s}{1}
	 \foreach \x in {#1}{
 	{\ifnum\value{r}=1
		\draw (0,\value{r}-1)--(\x,\value{r}-1); 
		\else
		{\ifnum\value{s}<\x
		\draw (\value{s},\value{r}-1)--(\x,\value{r}-1); 
		\fi}
		\fi}
	\draw (0,\value{r}) to (\x,\value{r});
   	\foreach \y in {0, ..., \x} {\draw (\y,\value{r})--(\y,\value{r}-1);}
	\addtocounter{r}{1}
	\setcounter{s}{\x}
 }}
 
\newcommand\Tableau[1]{
        \setcounter{r}{0}
        \setcounter{s}{0}
        \foreach \x [count = \c from 1] in {#1} {
		\foreach \y [count = \d from 1] in \x{
			\node at (\d-.5,\c-.5) {\scriptsize$\y$};
			\draw (\d,\c) to (\d,\c-1);
			{\ifnum\d=1
				\draw (0,\c) to (0,\c-1);
				\fi}
			\setcounter{r}{\d}
		}
		{\ifnum\c=1
			\draw (0,0)--(\value{r},0);
			\fi}
		\draw(0,\c) to (\value{r},\c);
		\setcounter{s}{\c}}}
\newcommand\smTableau[1]{
        \setcounter{r}{0}
        \setcounter{s}{0}
        \foreach \x [count = \c from 1] in {#1} {
		\foreach \y [count = \d from 1] in \x{
			\node at (\d-.5,\c-.5) {\small$\y$};
			\draw (\d,\c) to (\d,\c-1);
			{\ifnum\d=1
				\draw (0,\c) to (0,\c-1);
				\fi}
			\setcounter{r}{\d}
		}
		{\ifnum\c=1
			\draw (0,0)--(\value{r},0);
			\fi}
		\draw(0,\c) to (\value{r},\c);
		\setcounter{s}{\c}}}
\newcommand\nTableau[1]{
        \setcounter{r}{0}
        \setcounter{s}{0}
        \foreach \x [count = \c from 1] in {#1} {
		\foreach \y [count = \d from 1] in \x{
			\node at (\d-.5,\c-.5) {$\y$};
			\draw (\d,\c) to (\d,\c-1);
			{\ifnum\d=1
				\draw (0,\c) to (0,\c-1);
				\fi}
			\setcounter{r}{\d}
		}
		{\ifnum\c=1
			\draw (0,0)--(\value{r},0);
			\fi}
		\draw(0,\c) to (\value{r},\c);
		\setcounter{s}{\c}}}

\newcommand\sTableau[1]{
        \foreach \x [count = \c from 1] in {#1} {
		\foreach \y [count = \d from 1] in \x{
			\node at (\d-.5,\c-.5) {\tiny$\y$};
			\draw (\d,\c) to (\d,\c-1);
			{\ifnum\d=1
				\draw (0,\c) to (0,\c-1);
				\fi}
			\setcounter{r}{\d}
		}
		{\ifnum\c=1
			\draw (0,0)--(\value{r},0);
			\fi}
		\draw(0,\c) to (\value{r},\c);
		\setcounter{s}{\c}}}

\newcommand\XWord[2]{
 $\begin{tikzpicture}
		\foreach \x in {1, ..., #1}{\node (\x) at (-.35*\x, 0) {$s_\x$};}
		\node at (0, 0) {$\pi$};
		\foreach \x in {#2}{
		\node[inner sep=0pt, fill=black, opacity=.25] at (\x) {$\phantom{s_\x}$};
		\node[cross out, draw, inner sep=0pt] at (\x) {$\phantom{s_\x}$};
		}
		\end{tikzpicture}$}

\newcommand{\posleq}[1]{
	\hspace{0.1cm}
	\begin{tikzpicture}
	\draw (-0.8ex, -0.5ex) -- (0.8ex, -0.5ex);
	\draw (-0.8ex, 0.4ex) -- (0.7ex, -0.2ex);
	\draw (-0.8ex, 0.4ex) -- (0.7ex, 1ex);
	\draw (0.4ex,0.4ex) --(1.1ex, 0.4ex);
	\draw (0.75ex,0.75ex) --(0.75ex, 0.05ex);
	\end{tikzpicture}
	\hspace{0.1cm}
	}
\newcommand{\negleq}[1]{
	\hspace{0.1cm}
	\begin{tikzpicture}
	\draw (-0.8ex, -0.5ex) -- (0.8ex, -0.5ex);
	\draw (-0.8ex, 0.4ex) -- (0.7ex, -0.2ex);
	\draw (-0.8ex, 0.4ex) -- (0.7ex, 1ex);
	\draw (0.4ex,0.4ex) --(1.1ex, 0.4ex);
	\end{tikzpicture}
	\hspace{0.1cm}
	}
	
\newcommand{\zeroleq}[1]{
	\hspace{0.1cm}
	\begin{tikzpicture}
	\draw (-0.8ex, -0.5ex) -- (0.8ex, -0.5ex);
	\draw (-0.8ex, 0.4ex) -- (0.7ex, -0.2ex);
	\draw (-0.8ex, 0.4ex) -- (0.7ex, 1ex);
	\draw  (0.75ex,0.4ex) ellipse (0.2ex and 0.35ex);
	\end{tikzpicture}
	\hspace{0.1cm}
	}
	
\newcommand{\posgeq}[1]{
	\hspace{0.1cm}
	\begin{tikzpicture}
	\draw (-0.8ex, -0.5ex) -- (0.8ex, -0.5ex);
	\draw (0.8ex, 0.4ex) -- (-0.7ex, -0.2ex);
	\draw (0.8ex, 0.4ex) -- (-0.7ex, 1ex);
	\draw (-0.4ex,0.4ex) --(-1.1ex, 0.4ex);
	\draw (-0.75ex,0.75ex) --(-0.75ex, 0.05ex);
	\end{tikzpicture}
	\hspace{0.1cm}
	}
\newcommand{\neggeq}[1]{
	\hspace{0.1cm}
	\begin{tikzpicture}
	\draw (-0.8ex, -0.5ex) -- (0.8ex, -0.5ex);
	\draw (0.8ex, 0.4ex) -- (-0.7ex, -0.2ex);
	\draw (0.8ex, 0.4ex) -- (-0.7ex, 1ex);
	\draw (-0.4ex,0.4ex) --(-1.1ex, 0.4ex);
	\end{tikzpicture}
	\hspace{0.1cm}
	}
	
\newcommand{\zerogeq}[1]{
	\hspace{0.1cm}
	\begin{tikzpicture}
	\draw (-0.8ex, -0.5ex) -- (0.8ex, -0.5ex);
	\draw (0.8ex, 0.4ex) -- (-0.7ex, -0.2ex);
	\draw (0.8ex, 0.4ex) -- (-0.7ex, 1ex);
	\draw  (-0.75ex,0.4ex) ellipse (0.2ex and 0.35ex);
	\end{tikzpicture}
	\hspace{0.1cm}
	}

\newcommand{\posl}[1]{
	\hspace{0.1cm}
	\begin{tikzpicture}
	\draw (-0.8ex, 0.4ex) -- (0.7ex, -0.2ex);
	\draw (-0.8ex, 0.4ex) -- (0.7ex, 1ex);
	\draw (0.4ex,0.4ex) --(1.1ex, 0.4ex);
	\draw (0.75ex,0.75ex) --(0.75ex, 0.05ex);
	\end{tikzpicture}
	\hspace{0.1cm}
	}
\newcommand{\negl}[1]{
	\hspace{0.1cm}
	\begin{tikzpicture}
	\draw (-0.8ex, 0.4ex) -- (0.7ex, -0.2ex);
	\draw (-0.8ex, 0.4ex) -- (0.7ex, 1ex);
	\draw (0.4ex,0.4ex) --(1.1ex, 0.4ex);
	\end{tikzpicture}
	\hspace{0.1cm}
	}
	
\newcommand{\zerol}[1]{
	\hspace{0.1cm}
	\begin{tikzpicture}
	\draw (-0.8ex, 0.4ex) -- (0.7ex, -0.2ex);
	\draw (-0.8ex, 0.4ex) -- (0.7ex, 1ex);
	\draw  (0.75ex,0.4ex) ellipse (0.2ex and 0.35ex);
	\end{tikzpicture}
	\hspace{0.1cm}
	}
	
\newcommand{\posg}[1]{
	\hspace{0.1cm}
	\begin{tikzpicture}
	\draw (0.8ex, 0.4ex) -- (-0.7ex, 1ex);
	\draw (0.8ex, 0.4ex) -- (-0.7ex, -0.2ex);
	\draw (-0.4ex,0.4ex) --(-1.1ex, 0.4ex);
	\draw (-0.75ex,0.75ex) --(-0.75ex, 0.05ex);
	\end{tikzpicture}
	\hspace{0.1cm}
	}
\newcommand{\negg}[1]{
	\hspace{0.1cm}
	\begin{tikzpicture}
	\draw (0.8ex, 0.4ex) -- (-0.7ex, -0.2ex);
	\draw (0.8ex, 0.4ex) -- (-0.7ex, 1ex);
	\draw (-0.4ex,0.4ex) --(-1.1ex, 0.4ex);
	\end{tikzpicture}
	\hspace{0.1cm}
	}
	
\newcommand{\zerog}[1]{
	\hspace{0.1cm}
	\begin{tikzpicture}
	\draw (0.8ex, 0.4ex) -- (-0.7ex, -0.2ex);
	\draw (0.8ex, 0.4ex) -- (-0.7ex, 1ex);
	\draw  (-0.75ex,0.4ex) ellipse (0.2ex and 0.35ex);
	\end{tikzpicture}
	\hspace{0.1cm}
	}

\definecolor{plumb}{HTML}{8105C1}
\newcommand{\uoverline}[1]{\underline{#1}}
\newcommand{\Mk}[2]{\begin{tikzpicture}[xscale=.2]
        \setcounter{r}{1}
\foreach \x [count=\c from 1] in {#1}{\node (\c) at (\c,0) {\small $\x$};
        \setcounter{r}{\c}}
\node[below, outer sep=2pt, opacity=0] at (1) {\rule[1pt]{5pt}{1pt}};
\foreach \x  in {#2}{\node[below, plumb, outer sep=2pt] at (\x) {\rule[1pt]{5pt}{1pt}};}
\node[left] at (1,0) {\small $(\!$};\node[right] at (\value{r},0) {\small $\!)$};
\end{tikzpicture}}

\newcommand{\bred}[1]{\red{\bf #1}}

\newcommand{\rul}[1]{{\color{red}\underline{{#1}}}}

\makeatletter
\renewcommand{\@makefnmark}{\mbox{\textsuperscript{}}}
\makeatother

\title{Set-valued tableaux for Macdonald polynomials}
\author{
Zajj Daugherty\quad email:\ zdaugherty@gmail.com \\
Arun Ram\quad\quad\quad\ \ email:\ aram@unimelb.edu.au \\
\\
}
\date{}

\usetikzlibrary{arrows.meta}

\begin{document}

\maketitle

\vspace{-3em}
\begin{center}
{\sl In memory of Georgia Benkart}
\end{center}

\begin{abstract}
\noindent
Set-valued tableaux formulas play an important role in
Schubert calculus.
Using the box greedy reduced word for the construction of the Macdonald polynomials,
we convert the alcove walk formula for Macdonald
polynomials to a set-valued tableaux formula for Macdonald polynomials.
Our hope is that providing set-valued tableaux formulas for Macdonald 
polynomials will help
to strengthen the analogies and possible connections between the calculus
of Macdonald polynomials and Schubert calculus.
\end{abstract}


\section{Introduction}

The goal of this paper is to give set-valued tableaux formulas for Macdonald polynomials.
As evidenced by the work of \cite{Bu02}, \cite{KMY09}, \cite{MS19}, \cite{PS19}, \cite{Wei20}, \cite{Yu21} and others, 
set-valued tableaux formulas play an important role in
Schubert calculus and our hope is that providing set-valued tableaux formulas for Macdonald polynomials will help
to strengthen the analogies and possible connections between the calculus
of Macdonald polynomials and Schubert calculus.

Our set-valued tableaux formulas for Macdonald polynomials 
are derived by making a bijection between set-valued tableaux
with alcove walks.  Using this bijection and converting the alcove walks formula from \cite[Theorem 2.2]{RY08}
produces a set-valued tableaux formula for Macdonald polynomials.
We follow the framework of \cite[\S1 and \S2]{GR21}, which gives an exposition of the alcove walks formula
in the $GL_n$ context and analyzes the favourite reduced for the $n$-periodic permutation $u_\mu$
that is used to construct the relative Macdonald polynomial $E^z_\mu$.

Providing a set-valued tableaux formula for $E^z_\mu$ also provides a set-valued tableau formula for
the symmetric Macdonald polynomials $P_\lambda$ since if $\lambda$ is a partition then
(see \cite[(1.4),(1.5) and (1.10)]{GR21b})
$$P_\lambda = (const)\sum_{z\in S_n\lambda} E^z_\lambda.$$

The set-valued tableaux formulas given in Theorem \ref{pospowers} (respectively, Theorem \ref{negpowers}) are specializable
at $q=0$ and $t=0$ (respectively, $q=\infty$ and $t=\infty$).
The specializations of $E^z_\mu$ at $t=0$ and $t=\infty$ have interpretations in terms of crystals
for level 1 and level 0 affine Demazure characters (see, for example, \cite[Theorem C and Corollary D]{Kt16}
or \cite[(2.12), Theorem 1.1 and \S2.4]{MRY19}).  
It would be interesting to write the root operators for these crystals explicitly on the set-valued tableaux
(it is likely that the crystals of \cite[\S4]{Yu21} cover some cases of these crystals).

\section{Set-valued tableaux for Macdonald polynomials}

We begin with combinatorial definitions necessary for stating the set-valued tableaux formula
for Macdonald polynomials.  
To aid in processing these definitions, we encourage the reader to follow
Example \ref{runningexample} in tandem.

Fix $n\in \ZZ_{>0}$ and let $S_n$ denote the symmetric group on $n$ letters.
Identify $\mu=(\mu_1,\ldots, \mu_n)\in \ZZ_{\ge 0}^n$
with \emph{the set of boxes in $\mu$},
$$\mu = \{ (r,c)\ |\ \hbox{$r\in \{1, \ldots, n\}$ and $c\in \{1, \ldots, \mu_r\}$} \}.$$

\subsection{Definition of the shift and height statistics}

The minimal length permutation $v_\mu\in S_n$ such that $v_\mu \mu$ is weakly increasing
is given by
$$v_\mu(r) = 1+\#\{r'\in \{1, \ldots, r-1\}\ |\ \mu_{r'}\le \mu_r\}
+\#\{r'\in \{r+1, \ldots, n\}\ |\ \mu_{r'}<\mu_r\},$$
for $r\in \{1, \ldots, n\}$.
For a box $(r,c)\in \mu$ and $i\in \{1, \ldots, u_\mu(r,c)\}$ define
\begin{equation}
\mathrm{sh}(r,c) = \mu_r-c+1
\qquad\hbox{and}\qquad
\mathrm{ht}(i,r,c)
= v_\mu(r) - i.
\label{shhtdefn}
\end{equation}

\subsection{Definition of set-valued tableaux}

For $(r,c)\in \mu$ define
\begin{equation}
u_\mu(r,c) = \#\{ r'\in \{1, \ldots, r-1\}\ |\ \mu_{r'}<c\le \mu_r\} + \#\{r'\in \{r+1, \ldots, n\}\ |\ \mu_{r'}<c-1<\mu_r\}.
\label{eq:u(r,c)}
\end{equation}
The values $u_\mu(r,c)$ play a similar role to the flagging in the use of set-valued tableaux to give
formulas for Grothendieck polynomials (see \cite[\S2.6 and Theorem 7.1]{Wei20}).
A \emph{set-valued tableau $T$ of shape $\mu$} 
is a choice of subset $T(r,c)\subseteq \{1,\ldots, u_\mu(r,c)\}$ 
for each box $(r,c)\in \mu$.  More formally, a set-valued tableau $T$ of shape $\mu$
is a function 
$$T\colon \mu \to \{ \hbox{subsets of $\{1, \ldots, n\}$}\}
\qquad\hbox{such that}\quad
T(r,c) \subseteq \{1,\ldots, u_\mu(r,c)\}.$$
Let $\vert T\vert$ denote the total number of entries in $T$.

\subsection{Definition of $z^T_{(r,c)}$ and $x^T$}

Using one-line notation, let  $\gamma_n$ be the $n$-cycle in the symmetric group $S_n$ given by
$$\gamma_n 
= \begin{pmatrix} 
n &1 &2 &\cdots &n-1\end{pmatrix}.$$
For positive integers $k_1, \ldots, k_\ell$ such that $k_1+\cdots+k_\ell = n$, let 
$$\gamma_{k_1}\times \cdots \times \gamma_{k_\ell}
\qquad\hbox{be the disjoint
product of cycles in}\qquad
S_{k_1}\times \cdots \times S_{k_\ell}\subseteq S_n.
$$

Fix a permutation $z\in S_n$.  
Order the boxes of $\mu$ down columns and then left to right (i.e., in increasing order of the values $r+nc$).
Starting from the permutation $z$,
associate a permutation $z^T_{(r,c)}$ to each box $(r,c)$ of $\mu$ as follows: \hfil\break
If $(r,c)\in \mu$ 
and $T(r,c) =\{m_1, \ldots, m_p\}$ with $1\le m_1<\cdots <m_p\le u_\mu(r,c)$ 
define
\begin{equation}
z_{(r,c)}^T = z_{(r',c')}^T
(\gamma_{m_1-0}\times  \gamma_{m_2-m_1} \cdots \times \gamma_{m_p-m_{p-1}}
\times \gamma_{u_\mu(r,c)+1-m_p}\times (\gamma_1)^{\times (n-u_\mu(r,c)-1)}) \gamma_n^{-1},
\label{zTrcdefn}
\end{equation}
where $(r',c')\in \mu$ denotes the box before $(r,c)$ in $\mu$.
In the case that $(r,c)$ is the first box in $\mu$ we let $z^T_{(r',c')} = z$.
Then define 
\begin{equation}
x^T = \prod_{(r,c)\in \mu} x_{z^T_{(r,c)}(n)},
\label{xenddefn}
\end{equation}
where $z^T_{(r,c)}(n)$ is the value of the permutation $z^T_{(r,c)}$ at $n$.

The permutation $z^T_{\mathrm{init}}=z$ is the \emph{left key} (or \emph{initial direction}) and
the permutation $z^T_{\mathrm{fin}}$ corresponding to the last box in $\mu$ 
is the \emph{right key} (or \emph{final direction}).
These left and right keys play important roles in Schubert calculus: for example in 
\cite[\S3]{Yu21} and in the statement of Pieri-Chevalley formulas 
(compare \cite[Theorem 1]{FL94} and \cite[Corollary]{PR98}).

\subsection{Definition of the cov and maj statistics.}
Keep the notation $T(r,c) =\{m_1, \ldots, m_p\}$ with $m_1<\ldots<m_p$.
Define
\begin{align}
T^z_<(r,c) &= \{ m_k \in T(r,c)\ |\ 
 z^T_{(r,c)}(m_{k-1} ) < z^T_{(r,c)}(m_k)\}
\qquad\hbox{and} 
\nonumber \\
T^z_>(r,c) &= \{ m_k\in T(r,c)\ |\ 
z^T_{(r,c)}(m_{k-1} ) > z^T_{(r,c)}(m_k) \},
\label{pnfolds}
\end{align}
where we make the convention that $m_0=n$.
Then define
\begin{align*}
\mathrm{maj}_>^z(T) 
&= \sum_{(r,c)\in \mu} \mathrm{sh}(r,c)\cdot \vert T^z_>(r,c)\vert, \\
\mathrm{cov}_>^z(T) 
&= \Big( \sum_{(r,c)\in \mu} \sum_{m\in T^z_>(r,c)} \mathrm{ht}(m,r,c) \Big)
+ \hbox{$\frac12$} \big(\ \ell(z^T_{\mathrm{fin}}) - \ell(z^T_{\mathrm{init}}v_\mu^{-1}) - \vert T\vert \ \big), 
\end{align*}
where $\ell(w)$ denotes the length of the permutation $w$ in $S_n$.

\subsection{Set-valued tableaux formulas for $E^z_\mu$.}

The relative Macdonald polynomial $E^z_\mu = E^z_\mu(x_1, \ldots, x_n;q,t)$
is what is termed the ``permuted basement nonsymmetric Macdonald polynomial'' in \cite{Al16}
(see \cite[(3.7)]{GR21} for further references and a definition in terms of Cherednik-Dunkl operators
and the double affine Hecke algebra action on polynomials).

\begin{thm}  \label{pospowers}
Let $z\in S_n$ and $\mu\in \ZZ_{\ge 0}^n$.  Then the relative Macdonald polynomial 
$E^z_\mu$ is
$$E^z_\mu = \sum_T  q^{\mathrm{maj}_>^z(T)}t^{\mathrm{cov}_>^z(T)}
\Big(\prod_{(r,c)\in\mu} \prod_{m\in T(r,c)} \frac{1-t}{1-q^{\mathrm{sh}(r,c)}t^{\mathrm{ht}(m,r,c)}}\Big)
x^T,
$$
where the sum is over set-valued tableaux $T$ of shape $\mu$.
\end{thm}

Because the powers of $q$ and $t$ in formula in Theorem \ref{pospowers} are nonnegative integers, this formula
is well suited to the specializations $q=0$ and/or $t=0$.  The following Theorem provides an alternate
formula which is better for identifying the specializations of $E^z_\mu$ at $q=\infty$ and/or $t=\infty$ (i.e.,
specializations at $q^{-1}=0$ and $t^{-1}=0$).
\begin{align*}
\mathrm{maj}^z_<(T) 
&
= \sum_{(r,c)\in \mu} \mathrm{sh}(r,c)\cdot \vert T^z_<(r,c)\vert, \\
\mathrm{cov}^z_<(T) 
&= \Big( \sum_{(r,c)\in \mu} \sum_{m\in T^z_<(r,c)} \mathrm{ht}(m,r,c) \Big)
- \hbox{$\frac12$} \big(\ \ell(z^T_{\mathrm{fin}}) - \ell(z^T_{\mathrm{init}}v_\mu^{-1}) - \vert T\vert \ \big).
\end{align*}

\begin{thm}  \label{negpowers} 
Let $z\in S_n$ and $\mu\in \ZZ_{\ge 0}^n$.  Then the relative Macdonald polynomial 
$E^z_\mu$ is
$$E^z_\mu = \sum_T  q^{-\mathrm{maj}^z_<(T)}t^{-\mathrm{cov}^z_<(T)}
\Big(\prod_{(r,c)\in\mu} \prod_{m\in T(r,c)} 
\frac{1-t^{-1}}{1-q^{-\mathrm{sh}(r,c)}t^{-\mathrm{ht}(m,r,c)}}\Big)
x^T,
$$
where the sum is over set-valued tableaux $T$ of shape $\mu$.
\end{thm}

The proof of Theorem \ref{negpowers} is obtained from Theorem \ref{pospowers} by multiplying
numerator and denominator of each coefficient by
$$
\prod_{(r,c)\in \mu}\prod_{m\in T(r,c)} q^{-\mathrm{sh}(r,c)}t^{-\mathrm{ht}(m,r,c)}.
$$
The proof of Theorem \ref{pospowers} is obtained by considering the reduced word for the $n$-periodic
permutation
$$u_\mu = \prod_{(r,c)\in \mu} s_{u_\mu(r,c)-1}\cdots s_2s_1\pi
\qquad \hbox{(product taken in increasing order of boxes)}
$$
studied in \cite[Propositions 2.1 and 2.2]{GR21} (see also \cite[(2.4.3)]{Mac03}) and
using an easy bijection (see \S\ref{sec:conversion})
$$
\{\hbox{alcove walks of type $\vec u_\mu$}\} \longleftrightarrow \{\hbox{set-valued tableaux of shape $\mu$}\} 
$$
This bijection gives a conversion between the alcove walks formula for $E^z_\mu$
in \cite[Theorem 1.1(a)]{GR21} (see also \cite[Theorem 2.2]{RY08} for this alcove walks formula in a root system language) 
and the set-valued tableaux formula in Theorem \ref{pospowers}.

\begin{example} \label{runningexample}
Let $n=5 $ and $\mu = (0,4,5,1,4)$. Then the box arrangement associated to $\mu$ is as follows, with boxes placed on a labeled grid, filled with values $u_\mu(r,c)$, and marked in the lower corner by their cylindrical numbers $r+nc$.

\begin{equation}\label{eq:Box arrangement example}
\TikZ{[xscale=.75,yscale=-.75]	
	\foreach \x in {1,2,3,4,5}{
		\draw[black!40, very thin] (\x, -.2) to (\x, 5.2);
		\node[black!40,above] at (\x-.5, 0) {\tiny $\x$};}
	\draw[black!40, very thin] (0, -.2) to (0, 5.2);
	\foreach \x in {1,2,...,5}{
		\draw[black!40, thin] (-.2,\x) to (5.2,\x);
		\node[black!40,left] at (0, \x-.5) {\tiny $\x$};}
	\draw[black!40, thin] (-.2,0) to (5.2,0);
	\node[black!40, above] at (3-.5, -.5) {\scriptsize column $c$};
	\node[black!40, left ] at (-.5, 3-.5) {\scriptsize row $r$};
	\setcounter{bx}{1}
	\foreach \Col [count = \x from 1] in {{2,3,4,5}, {2,3,5},{2,3,5}, {2,3,5}, {3}}{
	\foreach \y in \Col {
		\setcounter{bx}{\x}
		\foreach \k in {1,2,3,4}{\addtocounter{bx}{\x}}
		\addtocounter{bx}{\y}
		\node[inner sep = 1pt, above right, dblue] at (\x-1, \y) {\tiny \thebx};
		\draw[dblue] (\x-1+.425, \y) .. controls +(0,-.325) .. +(-.425,-.325);
		}
		}
        \begin{scope}[thick]
        \smTableau{
        {},
        {{~1}, {~1}, {~2}, {~2}},
        {{~1},{~1},{~2},{~2},{~3}},
        {{~1}},
        {{~1},{~2},{~2},{~2}}}
	\draw[] (0,0) to (0,1) to (4,1) (4,2) to (5,2) (1,4) to (4,4);
	\end{scope}
}
\end{equation}
Box by box, the corresponding factors $\frac{(1-t)}{1-q^{\mathrm{sh}(r,c)}t^{\mathrm{ht}(i,r,c)}}$ are
$$
\def\arraystretch{1.1}
\frac{(1-t)}{1-q^{\mathrm{sh}}t^{\mathrm{ht}}}
=\begin{array}{|T|T|T|T|T|}
\multicolumn{1}{|T}{\TikZ{\draw[white] (0,0) to (0,.75);}}
\\\cline{1-4}
\begin{matrix} \frac{(1-t)}{(1-q^4t^2)}\end{matrix}
&\begin{matrix} \frac{(1-t)}{(1-q^3t^2)}\end{matrix} 
&\begin{matrix} \frac{(1-t)}{(1-q^2t^2)} \\ \frac{(1-t)}{(1-q^2t)} \end{matrix}
&\begin{matrix} \frac{(1-t)}{(1-qt^2)} \\ \frac{(1-t)}{(1-qt)} \end{matrix}
\\\cline{1-5}
\begin{matrix} \frac{(1-t)}{(1-q^5t^4)}\end{matrix} 
&\begin{matrix} \frac{(1-t)}{(1-q^4t^4)}\end{matrix}
&\begin{matrix} \frac{(1-t)}{(1-q^3t^4)} \\ \frac{(1-t)}{(1-q^3t^3)}\end{matrix}
&\begin{matrix} \frac{(1-t)}{(1-q^2t^4)} \\ \frac{(1-t)}{(1-q^2t^3)}\end{matrix}
&\begin{matrix} \frac{(1-t)}{(1-qt^4)} \\ \frac{(1-t)}{(1-qt^3)} \\ \frac{(1-t)}{(1-qt^2)} \end{matrix}
\\\cline{1-5}
\begin{matrix} \frac{(1-t)}{(1-qt)}\end{matrix}
\\\cline{1-4}
\begin{matrix} \frac{(1-t)}{(1-q^4t^3)}\end{matrix}
&\begin{matrix} \frac{(1-t)}{(1-q^3t^3)} \\ \frac{(1-t)}{(1-q^3t^2)}\end{matrix}
&\begin{matrix} \frac{(1-t)}{(1-q^2t^3)} \\ \frac{(1-t)}{(1-q^2t^2)} \end{matrix}
&\begin{matrix} \frac{(1-t)}{(1-qt^3)} \\ \frac{(1-t)}{(1-qt^2)} \end{matrix}\\\cline{1-4}
\end{array}
$$

\noindent 
One set-valued tableau of shape $\lambda$ is 
\begin{equation}\label{eq:favoriteT}
T =
\ \TikZ{[xscale=.85,yscale=-.55]
     \smTableau{
        {},
        {{\{1\}}, {\emptyset}, {\{1\}}, {\{2\}}},
        {{\{1\}},{\emptyset},{\{1,\! 2\}},{\emptyset},{\{1,\! 3\}}},
        {{\emptyset}},
        {{\emptyset},{\emptyset},{\{1,\! 2\}},{\emptyset}}}
	\draw[] (0,0) to (0,1) to (4,1) (4,2) to (5,2) (1,4) to (4,4);
}
=
\ \TikZ{[xscale=.55,yscale=-.55]
     \smTableau{
        {},
        {{1}, {}, {1}, {2}},
        {{1},{},{1\,\! 2},{},{1\,\! 3}},
        {{}},
        {{},{},{1\,\! 2},{}}}
	\draw[] (0,0) to (0,1) to (4,1) (4,2) to (5,2) (1,4) to (4,4);
}
\end{equation}
which has size 
$\vert T \vert = 1+1+1+2+2+1+2 = 10$ (we shall often omit the set notation as in the right-most picture above).
The product in Theorem \ref{pospowers} corresponding to this set-valued tableau $T$ is
$$\prod_{(r,c)\in \mu} \prod_{m\in T(r,c)} \frac{(1-t)}{1-q^{\mathrm{sh}(r,c)}t^{\mathrm{ht}(m,r,c)}}
=
\def\arraystretch{1.8}
\begin{array}{|c|c|c|c|c|}
\multicolumn{1}{|c}{~}
\\\cline{1-4}
\frac{(1-t)}{(1-q^4t^2)}
&1
&\frac{(1-t)}{(1-q^2t^2)}
&\frac{(1-t)}{(1-qt)}
\\\hline
\frac{(1-t)}{(1-q^5t^4)}
&1 &\frac{(1-t)}{(1-q^3t^4)}\frac{(1-t)}{(1-q^3t^3)}
&1 &\frac{(1-t)}{(1-qt^4)}\frac{(1-t)}{(1-qt^2)}
\\\hline
1 \\\cline{1-4}
1 &1 &\frac{(1-t)}{(1-q^2t^3)}\frac{(1-t)}{(1-q^2t^2)} &1\\\cline{1-4}
\end{array}
\ .
$$

Let $z = \id$ be the identity.  
The \emph{box-by-box permutation sequence of $T$ with initial direction $z$} 
(written in 1-line notation $w = (w(1)\ \cdots \ w(n))$) is
\begin{equation}\label{eq:zFavT}
z^T =
\ \TikZ{[xscale=1.35,yscale=-.6]
     \smTableau{
        {},
        {{\Mk{2,3,4,5,\bred{1}}{1}}, 	{\Mk{3,2,4,5,\bred{1}}{}}, 	{\Mk{1,4,2,5,\bred{3}}{1}}, 	{\Mk{2,3,1,4,\bred{5}}{2}}},
        {{\Mk{3,4,5,1,\bred{2}}{1}},		{\Mk{3,4,5,1,\bred{2}}{}},	{\Mk{4,2,5,3,\bred{1}}{1,2}},	{\Mk{2,4,3,5,\bred{1}}{}},	{\Mk{5,3,1,4,\bred{2}}{1,3}}
        },
        {{\Mk{3,5,1,2,\bred{4}}{}}},
        {{\Mk{3,1,2,4,\bred{5}}{}},{(3412\bred{5})},{\Mk{2,5,3,1,\bred{4}}{1,2}},{\Mk{2,3,5,1,\bred{4}}{}}}}
	\draw[] (0,0) to (0,1) to (4,1) (4,2) to (5,2) (1,4) to (4,4);
}.
\end{equation}
In each box we have underlined the positions specified by the entries of $T$; these positions are the
$m_1, \ldots, m_p$ in Section \ref{pnfolds}.  The red highlighted entry indicates $z^T_{(r,c)}(n)$ which is used
in the formula for $x^T$ in Section \ref{zTrcdefn}, and also as $z^T_{(r,c)}(m_0)$ in the definition of
$T^z_<(r,c)$ and $T^z_>(r,c)$ in Section \ref{pnfolds}.  So we have
\begin{equation*}
x^{T} =
\TikZ{[xscale=.6,yscale=-.6]
     \nTableau{
        {},
        {{x_{\bred{1}}}, {x_{\bred{1}}}, {x_{\bred{3}}}, {x_{\bred{5}}}},
        {{x_{\bred{2}}},{x_{\bred{2}}},{x_{\bred{1}}},{x_{\bred{1}}},{x_{\bred{2}}}},
        {{x_{\bred{4}}}},
        {{x_{\bred{5}}},{x_{\bred{5}}},{x_{\bred{4}}},{x_{\bred{4}}}}}
	\draw[] (0,0) to (0,1) to (4,1) (4,2) to (5,2) (1,4) to (4,4);
}
= x_1^4x_2^3x_3x_4^3x_5^3\ ;
\end{equation*}
and, following the definition in Section \ref{pnfolds} and
splitting $T$ into \emph{ascents} and \emph{descents} in the sequence of underlined numbers (preceded by the red number $z_{(r,c)}^T(n)$): 
\begin{equation*}
T^z_< = 
\ \TikZ{[xscale=.5,yscale=-.5]
     \nTableau{
        {},
        {{1}, {}, {}, {}},
        {{1},{},{},{},{1}},
        {{}},
        {{},{},{2},{}}}
	\draw[] (0,0) to (0,1) to (4,1) (4,2) to (5,2) (1,4) to (4,4);
}
\qquad\text{and}\qquad
T^z_> = \ \TikZ{[xscale=.5,yscale=-.5]
     \nTableau{
        {},
        {{}, {}, {1}, {2}},
        {{},{},{2},{},{3}},
        {{}},
        {{},{},{1},{}}}
	\draw[] (0,0) to (0,1) to (4,1) (4,2) to (5,2) (1,4) to (4,4);
}\ .
\end{equation*}
The initial and  final directions of $T$ are
\begin{equation*}
z^T_{\mathrm{init}} = \id = (12345)
\qquad\hbox{and}\qquad
z^T_{\mathrm{fin}} = z^T_{(3,5)} = (53142),
\end{equation*}
where $z^T_{(3,5)}$ indicates the permutation in box $(3,5)$ of $z^T$.  
Since $v_\mu = (13524)$ then
$$\ell(z^T_{\mathrm{fin}}) = 4+2+0+1 = 7
\quad \text{ and } \quad 
\ell(z^T_{\mathrm{init}}v^{-1}_\mu) = 3,$$
so that 
$$\hbox{$\frac12$}(\ell(z^T_{\mathrm{fin}})-\ell(z^T_{\mathrm{init}}v^{-1}_\mu)-\vert T\vert )
=\hbox{$\frac12$}(7-3-10) = -3.$$
Then
\begin{equation*}
q^{\mathrm{maj}^z_>(T)}t^{\mathrm{cov}^z_>(T)}
= t^{\frac12(7-3-10)}\cdot 
\def\arraystretch{1.1}
\begin{array}{|c|c|c|c|c|}
\multicolumn{1}{|c}{~}
\\\cline{1-4}
1 &1 &q^2t^2 &qt \\\hline
1 &1 &q^3t^3 &1 &qt^2 \\\hline
1 \\\cline{1-4}
1 &1 &q^2t^3 &1\\\cline{1-4}
\end{array}
=t^{-3} q^9 t^{11}
= q^9 t^8.
\end{equation*}
\end{example}

\section{Proof of Theorem \ref{pospowers}}

In this section we describe the conversion from alcove walks to set-valued tableaux.
We follow the framework of \cite[\S1]{GR21}, which gives an exposition of the alcove walks formula
in the $GL_n$ context and analyzes the favourite reduced for the $n$-periodic permutation $u_\mu$
that is used to construct the relative Macdonald polynomial $E^z_\mu$.

\subsection{Inversions and the box-greedy reduced word}

An \emph{$n$-periodic permutation} is a bijection $w\colon \ZZ\to \ZZ$ such that
$w(i+n) = w(i)+n$.  Since an $n$-periodic permutation $w$ is determined by the values
$w(1), \ldots, w(n)$ then
a permutation $w\in S_n$ extends uniquely to an $n$-periodic permutation.
For $i\in \{1, \ldots, n-1\}$, let $s_i$ be the transposition in $S_n$ switching $i$ and $i+1$.
Define $\pi\colon \ZZ\to \ZZ$ by $\pi(i) = i+1$.

An \emph{inversion} of an $n$-periodic permutation $w$ is an element of the set
$$\mathrm{Inv}(w) = \{ (i,k)\ |\ \hbox{$i\in \{1, \ldots, n\}$, $k\in \ZZ$, $i<k$ and $w(i)>w(k)$}\}.$$
If $i,j\in \{1, \ldots, n\}$ with $i<j$ and $\ell\in \ZZ$ then the 
\emph{shift} and \emph{height} of $(i,j+\ell n)$ are defined to be
$$\mathrm{sh}(i, j+\ell n) = \ell
\qquad\hbox{and}\qquad
\mathrm{ht}(i, j+\ell n) = j-i.
$$

Let $\mu = (\mu_1, \ldots, \mu_n)\in \ZZ^n_{\ge 0}$
and let $u_\mu$ be the $n$-periodic permutation defined
by
\begin{equation}
u_\mu(i) = v^{-1}_\mu(i) + n \mu_i , 
\qquad\hbox{for $i\in \{1, \ldots, n\}$.}
\end{equation}
Let $u_\mu(r,c)$ be as defined in \eqref{eq:u(r,c)}.
Following \cite[Prop. 2.2(1)]{GR21},
the \emph{box-greedy reduced word for $u_\mu$} (an abstract word in symbols $s_1,\ldots, s_{n-1}$ and $\pi)$ is
\begin{equation}
u^\square_\mu = \prod_{\mathrm{boxes}\ (r,c)\ \mathrm{in}\  dg(\mu)} (s_{u_\mu(r,c)}\cdots s_2s_1\pi).
\label{bgredwddefn}
\end{equation}
and the product is taken in increasing order of the values $r+nc$.  As in \cite[(2.12)]{GR21},
this reduced word induces a partition of the 
elements of $\mathrm{Inv}(u_\mu)$ according to the boxes $(r,c)$ of $\mu$, and an ordering of the
inversions in each box.
For $(r,c)\in \mu$ and $i\in \{1, \ldots, u_\mu(r,c)\}$,
$$\hbox{the $i$th inversion in the box $(r,c)$ is}\qquad
\beta(i,r,c) = (v_\mu(r), i+n\cdot(\mu_r-c+1)).$$
The inversion $\beta(i,r,c)$ has shift and height
\begin{equation}
\mathrm{sh}(r,c) = \mathrm{sh}(\beta(i,r,c)) = \mu_r-c+1
\qquad\hbox{and}\qquad
\mathrm{ht}(i,r,c) 
= \mathrm{ht}(\beta(i,r,c)) = v_\mu(r) - i.
\label{rttconversion}
\end{equation}

\begin{example} \label{boxgreedyexample}
If $\mu = (0,4,5,1,4)$ then $n=5$, 
the box-greedy reduced word for $u_\mu$ is
\def\YA{.5}\def\YB{1.5}\def\YC{3}\def\YD{4.5}\def\YE{6}\def\YF{7.9}
\begin{equation}
u^\square_\mu = (s_1\pi)^6(s_2s_1\pi)^7(s_3s_2s_1\pi) = 
\TikZ{[xscale=.75,yscale=-.65]	
\foreach \y in {2,3,4,5}{\node at (1,\y) { $s_1 \pi$};}
\foreach \y in {2,3}{\node at (2.25,\y) { $s_1 \pi$};}
\foreach \y in {5}{\node at (2.25,\y) { $s_2 s_1 \pi$};}
\foreach \y in {2,3,5}{\node at (3.75,\y) { $s_2 s_1 \pi$};}
\foreach \y in {2,3,5}{\node at (5.25,\y) { $s_2 s_1 \pi$};}
\foreach \y in {3}{\node at (6.95,\y) { $s_3 s_2 s_1 \pi$};}
\draw (\YA,.5) to (\YA,5.5);
\draw (\YB,1.5) to (\YB,5.5);
\draw (\YC,1.5) to (\YC, 3.5) (\YC, 4.5) to (\YC, 5.5);
\draw (\YD,1.5) to (\YD, 3.5) (\YD, 4.5) to (\YD, 5.5);
\draw (\YE,1.5) to (\YE, 3.5) (\YE, 4.5) to (\YE, 5.5);
\draw (\YF, 2.5) to (\YF, 3.5);
\foreach \x [count = \c from 1]  in {\YE, \YF, \YF, \YE, \YE}{
\draw (\YA, \c + .5) to (\x, \c + .5);}
}
\label{bgredwdexample}
\end{equation}
and
the inversion set of $u_\mu$ is 
$$
\def\arraystretch{1.1}
\mathrm{Inv}(u_\mu) = 
\text{\small $\begin{array}{|T|T|T|T|T|}
\multicolumn{1}{|T}{\TikZ{\draw[white] (0,0) to (0,.75);}}
\\\cline{1-4}
\begin{matrix} (3,1+4\cdot 5) \end{matrix}
&\begin{matrix} (3,1+3\cdot 5)  \end{matrix}
& \begin{matrix} (3,1+2\cdot 5) \\ (3,2+2\cdot 5) \end{matrix}
&\begin{matrix} (3, 1+1\cdot 5) \\ (3, 2+1\cdot 5)\end{matrix} 
\\\hline
\begin{matrix} (5,1+5\cdot 5) \end{matrix}
&\begin{matrix} (5,1+4\cdot 5) \end{matrix}\!
&\begin{matrix} (5,1+3\cdot 5) \\ (5, 2 + 3\cdot 5) \end{matrix}
&\begin{matrix} (5, 1+2\cdot 5) \\ (5, 2+2\cdot 5)  \end{matrix}
&\begin{matrix} (5,1+1\cdot 5) \\ (5, 2+1\cdot 5) \\ (5, 3+1\cdot 5) \end{matrix}
\\\hline
\begin{matrix} (2, 1+1\cdot 5)\end{matrix} &\multicolumn{1}{T}{\TikZ{\draw[white] (0,0) to (0,1);}}
\\\cline{1-4}
\begin{matrix} (4, 1+4\cdot 5)\end{matrix}
&\begin{matrix} (4,1+3\cdot 5) \\ (4, 2+3\cdot 5) \end{matrix} 
&\begin{matrix} (4, 1+2\cdot 5) \\ (4, 2 + 2\cdot 5) \end{matrix} 
&\begin{matrix} (4, 1+1\cdot 5) \\ (4, 2+ 1\cdot 5) \end{matrix} \\
\cline{1-4}
\end{array}$}
$$
\end{example}

\subsection{Alcove walks and permutation sequences}

Let $z\in S_n$ and 
choose a subset $F$ of the $s_i$ factors in $u^\square_\mu = w_1w_2\cdots w_r$ to cross out.  
The corresponding \emph{alcove walk of type $(z, u_\mu^\square)$}
is the sequence $p(F)=(p_0, p_1, \ldots, p_r)$ of $n$-periodic permutations 
given by
$$\hbox{$p_0 = z$}\quad\hbox{and}\quad
p_k = \begin{cases}
p_{k-1}\pi, &\hbox{if $w_k = \pi$}, \\
p_{k-1}, &\hbox{if $w_k \in F$ (so that the factor $w_k$ is crossed out)}, \\
p_{k-1}w_k, &\hbox{if $w_k \not\in F$ (so that the factor $w_k$ is not crossed out)}.
\end{cases}
$$
Following \cite[(1.14)]{GR21},
the \emph{permutation sequence of $p(F)$} is the sequence of permutations in $S_n$ given by
$$\overline{p}(F) = (\overline{p_0}, \overline{p_1}, \ldots, \overline{p_r}),
\qquad\hbox{where}\quad
\overline{w}(i) = w(i) \bmod n.
$$
Note that $\overline{w_1w_2} = \overline{w_1}\,\overline{w_2}$ and 
$\overline{\pi} = c_n = s_1\cdots s_{n-1}$ and
$\overline{s_j} = s_j$ for $j\in \{1, \ldots, n-1\}$.

\subsection{Converting to set-valued tableaux} \label{sec:conversion}

The bijection between alcove walks and set-valued tableaux is
$$
\begin{matrix}
\{\hbox{alcove walks of type $u^\square_\mu$}\} &\longleftrightarrow &\{\hbox{set-valued tableaux of shape $\mu$}\} \\
p(F) &\longmapsto &T
\end{matrix}
$$
where, in the set-valued tableau $T$ corresponding to the alcove walk $p(F)$,
the set $T(r,c)$ in box $(r,c)$ specifies which factors are crossed out in that box: if a box in $T$ contains $i$ then 
delete $s_i$ in the corresponding word. 

\begin{example} \label{pfbyexample}
In the case $\mu = (0,4,5,1,4)$, where the box-greedy reduced word for $u_\mu$ is as given in 
\eqref{bgredwdexample},
there are $5+2\cdot 8+3=24$ factors of the form $s_i$ in $u^\square_\mu$ and 
so there are a total of $2^{24}$ alcove walks of type $(z,u^\square_\mu)$ 
(for any fixed permutation $z\in S_n$).
Each alcove walk corresponds to a choice of $s_i$ factors in $u^\square_\mu$ to cross out.  

The tableau $T$ in Example \ref{runningexample} corresponds to the subset 
\def\YA{.5}\def\YB{1.5}\def\YC{3}\def\YD{4.5}\def\YE{6}\def\YF{8}
\begin{equation}
F = 
\TikZ{[xscale=.8,yscale=-.7]	
\foreach \y in {2,3}{\node at (1,\y) {\XWord{1}{1}};}
\foreach \y in {4,5}{\node at (1,\y) { $s_1\pi$};}
\foreach \y in {2,3}{\node at (2.25,\y) { $s_1 \pi$};}
\foreach \y in {5}{\node at (2.25,\y) { $s_2 s_1 \pi$};}
\foreach \y in {2}{\node at (3.75,\y) {\XWord{2}{1}};}
\foreach \y in {3,5}{\node at (3.75,\y) {\XWord{2}{1,2}};}
\foreach \y in {2}{\node at (5.25,\y) {\XWord{2}{2}};}
\foreach \y in {3,5}{\node at (5.25,\y) { $s_2 s_1 \pi$};}
\foreach \y in {3}{\node at (7,\y) {\XWord{3}{1,3}};}
\draw (\YA,.5) to (\YA,5.5);
\draw (\YB,1.5) to (\YB,5.5);
\draw (\YC,1.5) to (\YC, 3.5) (\YC, 4.5) to (\YC, 5.5);
\draw (\YD,1.5) to (\YD, 3.5) (\YD, 4.5) to (\YD, 5.5);
\draw (\YE,1.5) to (\YE, 3.5) (\YE, 4.5) to (\YE, 5.5);
\draw (\YF, 2.5) to (\YF, 3.5);
\foreach \x [count = \c from 1]  in {\YE, \YF, \YF, \YE, \YE}{
\draw (\YA, \c + .5) to (\x, \c + .5);}
}\ ,
\label{alcwalkexample}
\end{equation}
which has alcove walk
$$p(F)=(p_0, p_1, \ldots, p_{37}) = (z, z, z\pi, z\pi, z\pi^2, z\pi^2 s_1, z\pi^2 s_1 \pi, z\pi^2 s_1 \pi s_1, z\pi^2 s_1 \pi s_1 \pi,   \ldots)$$
(there is a repeat entry in $p(F)$ each time there is an $s_i$ crossed out in $F$).
Using one line notation $w=(w(1)\ w(2)\ \cdots w(n))$ for $n$-periodic permutations (underlining 2-digit terms for emphasis), a box-by-box formulation in the case that $z = \id = (12345)$ is
{\def\arraystretch{1.2}
\setlength{\arraycolsep}{3pt}
\begin{equation*}
p(F) = \quad \text{\small$
\begin{array}{|T|T|T|T|T|}
\multicolumn{1}{|T}{\TikZ{\draw[white] (0,0) to (0,.75);}}
\\\cline{1-4}
\begin{matrix} (1\ 2\ 3\ 4\ 5) \\ (2\ 3\ 4\ 5\ 6)  \end{matrix} 
&\begin{matrix}  (6\ 3\ 7\ 9\ \underline{10}) \\ (3\ 7\ 9\ \underline{10}\ \uoverline{11}) \end{matrix}
&\begin{matrix} 
(3\ \underline{11}\ 9\ \uoverline{12}\ \underline{15}) \\
(3\ \underline{11}\ 9\ \uoverline{12}\ \underline{15}) \\
(\underline{11}\ 9\ \uoverline{12}\ \underline{15}\ 8) 
\end{matrix}
&\begin{matrix}
(\uoverline{12}\ \underline{15}\ 8\ \uoverline{16}\ \underline{14}) \\
(\uoverline{15}\ \underline{12}\ 8\ \uoverline{16}\ \underline{14}) \\
(\underline{12}\ 8\ \uoverline{16}\ \underline{14}\ \uoverline{20})
\end{matrix}
\\\hline
\begin{matrix} (2\ 3\ 4\ 5\ 6) \\ (3\ 4\ 5\ 6\ 7) \end{matrix} 
&\begin{matrix}   (7\ 3\ 9\ \underline{10}\ \uoverline{11}) \\ 
(3\ 9\ \underline{10}\ \uoverline{11}\ \underline{12})  \end{matrix}
&\begin{matrix} 
(\underline{11}\ 9\ \uoverline{12}\ \underline{15}\ 8) \\
(\underline{11}\ 9\ \uoverline{12}\ \underline{15}\ 8) \\
(9\ \uoverline{12}\ \underline{15}\ 8\ \underline{16}) 
\end{matrix}
&\begin{matrix} 
(\underline{12}\ \uoverline{16}\ 8\ \underline{14}\ \uoverline{20}) \\
(\underline{16}\ \uoverline{12}\ 8\ \underline{14}\ \uoverline{20}) \\
(\uoverline{12}8\underline{14}\ \uoverline{20}\ \underline{21})
\end{matrix}
&\begin{matrix} 
(\underline{12}\ 8\ \uoverline{20}\ \underline{21}\ \uoverline{19}) \\
(\underline{12}\ \uoverline{20}\ 8\ \underline{21}\ \uoverline{19}) \\
(\underline{12}\ \uoverline{20}\ 8\ \underline{21}\ \uoverline{19}) \\
(\uoverline{20}\ 8\ \underline{21}\ \uoverline{19}\ \underline{17})
\end{matrix} 
\\\hline
\begin{matrix} (4\ 3\ 5\ 6\ 7) \\ (3\ 5\ 6\ 7\ 9)\end{matrix} 
\\\cline{1-4}
\begin{matrix} (5\ 3\ 6\ 7\ 9) \\ (3\ 6\ 7\ 9\ \underline{10}) \end{matrix}
&\begin{matrix} (3\ \underline{10}\ 9\ \uoverline{11}\ \underline{12}) \\
(\underline{10}\ 3\ 9\ \uoverline{11}\ \underline{12}) \\ 
(3\ 9\ \underline{11}\ \uoverline{12}\ \underline{15})
\end{matrix}
&\begin{matrix} 
(9\ \uoverline{12}\ \underline{15}\ 8\ \uoverline{16}) \\
(9\ \uoverline{12}\ \underline{15}\ 8\ \uoverline{16}) \\
(\uoverline{12}\ \underline{15}\ 8\ \uoverline{16}\ \underline{14}) 
\end{matrix}
&\begin{matrix} 
(\uoverline{12}\ \underline{14}\ 8\ \uoverline{20}\ \underline{21}) \\
(\uoverline{14}\ \underline{12}\ 8\ \uoverline{20}\ \underline{21}) \\
(\underline{12}\ 8\ \uoverline{20}\ \underline{21}\ \uoverline{19})
\end{matrix}\\\cline{1-4}
\end{array}
$}
\end{equation*}}
The permutation sequence is
$$\overline{p}(F)=(\overline{p_0}, \overline{p_1}, \ldots, \overline{p_{37}}) 
	= (z, z, zc_n, zc_n, zc_n^2, zc_n^2 s_1, zc_n^2 s_1 c_n, zc_n^2 s_1 c_n s_1, zc_n^2 s_1 c_n s_1 c_n,   \ldots)$$
(obtained by replacing each $\pi$ by the $n$-cycle $c_n = s_1\cdots s_{n-1} = \gamma_n^{-1}$).
Namely, $\overline{p}(F)$ is obtained by taking all the values in $p(F)$ mod $n$ (here $n=5$ for $\mu = (0,4,1,5,4)$).
\begin{equation}
\overline{p}(F) = \quad 
\text{\small
\setlength{\arraycolsep}{3pt}
$
\begin{array}{|T|T|T|T|T|}
\multicolumn{1}{|T}{\TikZ{\draw[white] (0,0) to (0,.75);}}
\\\cline{1-4}
\begin{matrix} (\blue{\bf 12}345) \\ (23451)  \end{matrix} 
&\begin{matrix}  (13245) \\ (32451) \end{matrix}
&\begin{matrix} (31425) \\
(\blue{\bf 31}425) \\
(14253) 
\end{matrix}
&\begin{matrix}
(2\blue{\bf 53}14)  \\
(52314)  \\
(23145)  
\end{matrix}
\\\hline
\begin{matrix} (\blue{\bf 23}451) \\ (34512)  \end{matrix} 
&\begin{matrix}   (23451) \\ 
(34512)  \end{matrix}
&\begin{matrix} 
(1\blue{\bf 42}53)  \\
(\blue{\bf 14}253)  \\
(42531)  
\end{matrix}
&\begin{matrix} 
(21345)  \\
(12345)  \\
(24351) 
\end{matrix}
&\begin{matrix} 
(23\blue{\bf 51}4)  \\
(25314)  \\
(\blue{\bf 25}314)  \\
(53142)  \\
\end{matrix} 
\\\hline
\begin{matrix} (43512) \\ (35124) \end{matrix} 
\\
\cline{1-4}
\begin{matrix} (53124) \\ (31245) \end{matrix}
&\begin{matrix} (35412) \\
(53412) \\ 
(34125)
\end{matrix}
&\begin{matrix} 
(4\blue{\bf 25}31)  \\
(\blue{\bf 42}531)  \\
(25314)  
\end{matrix}
&\begin{matrix} 
(24351)  \\
(42351)  \\
(23514)  
\end{matrix}\\\cline{1-4}
\end{array}
$}
\label{alcwalkexample}
\end{equation}
The \emph{box-by-box permutation sequence} 
of $F$ is obtained by recording the last permutation in each box, which is $z^T$ in \eqref{eq:zFavT}. 
The elements $z^T_{(r,c)}$ in each box of $z^T$ are also obtained from $z=z^T_{\mathrm{init}}$ and $T$ by
the formula in Section \ref{zTrcdefn}. 

We have used the blue highlighted numbers in \eqref{alcwalkexample}
to record the positions $i$, $i+1$ when there is a repeat entry
in the sequence $\overline{p}(F)$ coming from a crossed out $s_i$.  
These record the information of $F$.  This same information is translated into the underlines in $z^T$ 
in \eqref{eq:zFavT}, and the underlines exactly specify the entries of $T$.
The set-valued tableau indicating the positions of the underlined numbers in each box is the tableau in \eqref{eq:favoriteT}.
\end{example}

\subsection{Conversion of the statistics} 

The proof of Theorem \ref{pospowers} is completed by matching up the statistics in Theorem \ref{pospowers}
with the statistics that appear in \cite[Theorem 3.1]{RY08}.  In Example \ref{pfbyexample}, we have used the 
blue highlighted numbers in \eqref{alcwalkexample} and the underlines in \eqref{eq:zFavT} 
to illustrate how the information of 
the `folds' in \cite[Theorem 3.1]{RY08} is equivalent to the information of the set valued tableau $T$.
In the context of \cite[(2.36)]{RY08}, the set $T(r,c)$ exactly records the folds coming from the box $(r,c)$.
By \eqref{rttconversion},
the factor
$\prod_{(r,c)\in \mu} \prod_{m\in T(r,c)} \frac{(1-t)}{1-q^{\mathrm{sh}(r,c)}t^{\mathrm{ht}(m,r,c)}}$
which appears in the product in Theorem \ref{pospowers} corresponds to
$$
\hbox{the factors}\quad
\left(\prod_{k\in f^+(p)} \frac{(1-t_{\beta^\vee_k})}{1-q^{\langle -\beta^\vee_k, \rho_c\rangle}} \right)
\left(\prod_{k\in f^-(p)} \frac{(1-t_{\beta^\vee_k})}{1-q^{\langle -\beta^\vee_k, \rho_c\rangle}} \right)
\quad
\hbox{in \ \cite[Theorem 3.1]{RY08}.}
$$
The sets $T^z_{<}$ and $T^z_{>}$ correspond to the sets $f^+(p)$ and $f^-(p)$ of positive and negative folds in \cite[Theorem 3.1]{RY08}.
The permutation $z^T_{\mathrm{fin}}$ is the permutation denoted $\varphi(p)$ in \cite[Theorem 3.1]{RY08}
(denoted $z_r$ in \cite[(1.14)]{GR21}),
and the permutation $z^T_{\mathrm{init}}v_\mu^{-1}$ is the permutation $m$ of \cite[Theorem 3.1]{RY08}.
With these conversions, the last term 
$\hbox{$\frac12$} \big(\ \ell(z^T_{\mathrm{fin}}) - \ell(z^T_{\mathrm{init}}v_\mu^{-1}) - \vert T\vert \ \big)$
in the definition of $\mathrm{cov}_>^z(T)$ in Section \ref{pnfolds}
corresponds to the combination of the 
factors $t^{\frac12}_{\varphi(p)}$ and $t^{-\frac12}_{\beta^\vee_k}$  in \cite[Theorem 3.1]{RY08} 
and the factor $t^{\frac12}_m$ in \cite[Remark 3.2]{RY08}.

\begin{example}
Suppose that
$n=14$, $\mu = (\mu_1, \ldots, \mu_{14})\in (\ZZ_{\ge 0})^{14}$, $(r,c)$ is a box in $\mu$
and $(r',c')$ is the box before $(r,c)$ in $\mu$.
Say 
$$u_{(r,c)}=11,
\qquad
T(r,c) = \{3,4,8,10\}
\qquad\hbox{and}\qquad
z^T_{(r',c')} = (u_1 u_2 u_3 u_4 u_5 u_6 u_7 u_8 u_9 u_{10} u_{11} u_{12} u_{13} u_{14}),
$$  
so that $u_i = z^T_{(r',c')}(i)$ for $i\in \{1, \ldots, 14\}$.  
Then the box $(r,c)$ of $F$ contains
\begin{align*}
s_{11}&\cancel{s_{10}}s_9\cancel{s_8}s_7s_6s_5\cancel{s_4}\cancel{s_3}s_2s_1\pi
= (\gamma_3\times \gamma_1\times \gamma_4\times \gamma_2\times \gamma_2\times \gamma_1\times \gamma_1)\pi
\\
&= (\gamma_{3-0}\times \gamma_{4-3}\times \gamma_{8-4}\times 
\gamma_{10-8}\times \gamma_{12-10}\times(\gamma_1)^{\times(14-12)})\pi
\end{align*}
and
the entries in box $(r,c)$ of the permutation sequence $\overline{p}(F)$ are
\begin{align*}
z^T_{(r',c')}s_{11} = z^T_{(r',c')}s_{11}\cancel{s_{10}} &= (u_1 u_2 u_3 u_4 u_5 u_6 u_7 u_8 u_9 u_{10} u_{12} u_{11} u_{13} u_{14}) \\
z^T_{(r',c')}s_{11}\cancel{s_{10}}s_9 = z^T_{(r',c')}s_{11}\cancel{s_{10}}s_9\cancel{s_8} &= (u_1 u_2 u_3 u_4 u_5 u_6 u_7 u_8  u_{10} u_9 u_{12} u_{11} u_{13} u_{14}) \\
z^T_{(r',c')}s_{11}\cancel{s_{10}}s_9\cancel{s_8} s_7
&= (u_1 u_2 u_3 u_4 u_5 u_6 u_8 u_7 u_{10} u_9 u_{12} u_{11} u_{13} u_{14}) \\
z^T_{(r',c')}s_{11}\cancel{s_{10}}s_9\cancel{s_8} s_7s_6
&= (u_1 u_2 u_3 u_4 u_5 u_8 u_6  u_7 u_{10} u_9 u_{12} u_{11} u_{13} u_{14}) \\
z^T_{(r',c')}s_{11}\cancel{s_{10}}s_9\cancel{s_8} s_7s_6s_5
&= (u_1 u_2 u_3 u_4 u_8 u_5  u_6  u_7 u_{10} u_9 u_{12} u_{11} u_{13} u_{14}) \\
z^T_{(r',c')}s_{11}\cancel{s_{10}}s_9\cancel{s_8} s_7s_6s_5\cancel{s_4}\cancel{s_3}s_2
&= (u_1 u_3 u_2  u_4 u_8 u_5  u_6  u_7 u_{10} u_9 u_{12} u_{11} u_{13} u_{14}) \\
z^T_{(r',c')}s_{11}\cancel{s_{10}}s_9\cancel{s_8} s_7s_6s_5\cancel{s_4}\cancel{s_3}s_2 s_1
&= ( u_3 u_1 u_2  u_4 u_8 u_5  u_6  u_7 u_{10} u_9 u_{12} u_{11} u_{13} u_{14}) \\
z^T_{(r,c)}
&= ( u_1 u_2  \rul{u_4} \rul{u_8} u_5  u_6  u_7 \rul{u_{10}} u_9 \rul{u_{12}} u_{11} u_{13} u_{14} \bred{u_3} )
\end{align*}
and
\begin{align*}
z^T_{(r,c)} 
&= z^T_{(r',c')} s_{11}\cancel{s_{10}}s_9\cancel{s_8}s_7s_6s_5\cancel{s_4}\cancel{s_3}s_2s_1 \gamma_{14}^{-1} \\
&= z^T_{(r',c')} (\gamma_3\times \gamma_1\times \gamma_4\times \gamma_2\times \gamma_2\times \gamma_1\times \gamma_1) \gamma_{14}^{-1} \\
&= z^T_{(r',c')} (\gamma_{3-0}\times \gamma_{4-3}\times \gamma_{8-4}\times 
\gamma_{10-8}\times \gamma_{12-10}\times
(\gamma_1)^{\times(14-12)}) \gamma_{14}^{-1}
\end{align*}
so that
$z^T_{(r,c)}(n) = u_3$,
$z^T_{(r,c)}(3) = u_4$,
$z^T_{(r,c)}(4) = u_8$,
$z^T_{(r,c)}(8) = u_{10}$,
$z^T_{(r,c)}(10) = u_{12}$.
Then $T^z_<(r,c)$ and $T^z_>(r,c)$ are determined by
$$\begin{array}{llll}
\hbox{if $u_3<u_4$\quad then}\quad &3\in T^z_<(r,c)\quad
&\hbox{and otherwise}\quad &3\in T^z_>(r,c); \\
\hbox{if $u_4<u_8$\quad then}\quad &4\in T^z_<(r,c)
&\hbox{and otherwise}\quad &4\in T^z_>(r,c); \\
\hbox{if $u_8<u_{10}$\quad then}\quad &8\in T^z_<(r,c)
&\hbox{and otherwise}\quad &8\in T^z_>(r,c); \\
\hbox{if $u_{10}<u_{12}$\quad then}\quad &10\in T^z_<(r,c)
&\hbox{and otherwise}\quad &10\in T^z_>(r,c).
\end{array}
$$
\end{example}

\section{Example: Computing $E_{(2,2,1,1,0,0)}$}

Here we compute $E_{(2,2,1,1,0,0)} = E^{\id}_{(2,2,1,1,0,0)}$ using the set-valued tableaux formula
of Theorem \ref{poposwers}.  
For comparison, this same Macdonald polynomial is computed in terms of nonattacking fillings
in \cite[example following Theorem 1.1]{GR21} and in terms of multiline queues in \cite[Figure 4]{CMW18}.

Let
$$
z = \id = (123456),
\quad\hbox{and}\quad
\mu = (2,2,1,1,0,0)
\quad\hbox{with}\quad
v_\mu = (563412)
\quad\hbox{and}\quad 
\ell(v_\mu) = 12
.$$
Then
$$u_\mu = \pi^4(s_2s_1\pi)^2 = 
\TikZ{[xscale=.5, yscale=-.5]
\draw (0,0) to (0,6) (1,0) to (1,4) (3.5,0) to (3.5,2); 
\foreach \x [count=\c from 0] in {3.5,3.5,3.5,1,1}{
	\draw (0,\c) to (\x,\c);
	}
\foreach \y in {1,2,3,4}{\node at (.5, \y-.5) {\footnotesize$\pi$};}
\node at (2.25, 1-.5) {$s_2 s_1 \pi$};
\node at (2.25, 2-.5) {$s_2 s_1 \pi$};}
\quad\hbox{and}\quad
\mathrm{Inv}(u_\mu)
= 
\TikZ{[xscale=.65, yscale=-.65]
\draw (0,0) to (0,5.25) (0,5.5) to (0,5.75) (1,0) to (1,5) (4.3,0) to (4.3,3); 
\draw[densely dotted] (0,5.25) to (0,5.5);
\foreach \y in {.75,2.25,3.5,4.5}{\node at (.5, \y) {$\emptyset$};}
\foreach \x/\y in {4.3/0, 4.3/1.5, 4.3/3, 1/4, 1/5}{\draw (0,\y) to (\x,\y);}
\node at (2.65, .75) {$\displaystyle 
	\begin{matrix} (5,1+1\cdot 6) \\ (5, 2+1\cdot 6) \end{matrix}$};
\node at (2.65, 2.25) {$\displaystyle 
	\begin{matrix} (6,1+1\cdot 6) \\ (6, 2+1\cdot 6) \end{matrix}$};
}
$$
Box by box, the corresponding factors $\frac{(1-t)}{1-q^{\mathrm{sh}(r,c)}t^{\mathrm{ht}(i,r,c)}}$ are
$$
\def\arraystretch{1.1}
\frac{(1-t)}{1-q^{\mathrm{sh}}t^{\mathrm{ht}}}
=\begin{array}{|T|T|T|T|T|}
\cline{1-2}
\begin{matrix} \phantom{\frac{(1-t)}{(1-q^4t^2)} } \end{matrix}
&\begin{matrix} \frac{(1-t)}{(1-qt^4)} \\ \frac{(1-t)}{(1-qt^3)} \end{matrix}
\\\cline{1-2}
\begin{matrix} \phantom{\frac{(1-t)}{(1-q^4t^2)} } \end{matrix}
&\begin{matrix} \frac{(1-t)}{(1-qt^4)} \\ \frac{(1-t)}{(1-qt^3)} \end{matrix}
\\\cline{1-2}
\begin{matrix} \phantom{\frac{(1-t)}{(1-qt)}} \end{matrix}
\\\cline{1-1}
\begin{matrix} \phantom{\frac{(1-t)}{(1-qt)}} \end{matrix}
\\\cline{1-1}
\end{array}
$$
For a set valued tableau $T$ let
$$
\wt(T) = q^{\mathrm{maj}_>^z(T)}t^{\mathrm{cov}_>^z(T)}
\Big(\prod_{(r,c)\in\mu} \prod_{m\in T(r,c)} \frac{1-t}{1-q^{\mathrm{sh}(r,c)}t^{\mathrm{ht}(m,r,c)}}\Big)
$$
where
$$\mathrm{cov}_>^z(T) 
= \Big( \sum_{(r,c)\in \mu} \sum_{m\in T^z_>(r,c)} \mathrm{ht}(m,r,c) \Big) 
+ \hbox{$\frac12$} \big(\ \ell(z^T_{\mathrm{fin}}) - \ell(z^T_{\mathrm{init}}v_\mu^{-1}) - \vert T\vert \ \big).
$$
For this example, $\mathrm{sh}(1,2) = \mathrm{sh}(2,2) = 1$ so that
$\mathrm{maj}_>^z(T) = \vert T^z_>(r,c)\vert$.

There are 16 terms in the expansion of $E_{(2,2,1,1,0,0)}$ corresponding to the following permutation sequences
$z^T$.

$$
z^{(1)} =
\ \TikZ{[xscale=1.55,yscale=-.6]
     \smTableau{
        {{\Mk{2,3,4,5,6,\bred{1}}{}}, 	{\Mk{6,1,2,3,4,\bred{5}}{1,2}}
        },
        {{\Mk{3,4,5,6,1,\bred{2}}{}},	{\Mk{1,2,3,4,5,\bred{6}}{1,2}}
        },
        {{\Mk{4,5,6,1,2,\bred{3}}{}}
        },
        {{\Mk{5,6,1,2,3,\bred{4}}{}}
        }
        }
	\draw[] (0,0) to (0,6);
}
\qquad\hbox{with}\qquad
\begin{array}{l}
x^T = x_1x_2x_3x_4x_5x_6 \\
\wt(T) = q^2\frac{(1-t)}{(1-qt^3)}\frac{(1-t)}{(1-qt^4)}\frac{(1-t)}{(1-qt^4)}\frac{(1-t)}{(1-qt^5)} \\
\mathrm{cov}^z_>(T) = 3+5+\frac12(0-12-4)=0.
\end{array}
$$
$$
z^{(2)} =
\ \TikZ{[xscale=1.55,yscale=-.6]
     \smTableau{
        {{\Mk{2,3,4,5,6,\bred{1}}{}}, 	{\Mk{6,1,2,3,4,\bred{5}}{1,2}}
        },
        {{\Mk{3,4,5,6,1,\bred{2}}{}},	{\Mk{2,1,3,4,5,\bred{6}}{1}}
        },
        {{\Mk{4,5,6,1,2,\bred{3}}{}}
        },
        {{\Mk{5,6,1,2,3,\bred{4}}{}}
        }
        }
	\draw[] (0,0) to (0,6);
}
\qquad\hbox{with}\qquad
\begin{array}{l}
x^T = x_1x_2x_3x_4x_5x_6 \\
\wt(T) = q^2 t \frac{(1-t)}{(1-qt^3)}\frac{(1-t)}{(1-qt^4)}\frac{(1-t)}{(1-qt^5)} \\
\mathrm{cov}^z_>(T) = 3+5+\frac12(1-12-3)=1.
\end{array}
$$
$$
z^{(3)} =
\ \TikZ{[xscale=1.55,yscale=-.6]
     \smTableau{
        {{\Mk{2,3,4,5,6,\bred{1}}{}}, 	{\Mk{6,1,2,3,4,\bred{5}}{1,2}}
        },
        {{\Mk{3,4,5,6,1,\bred{2}}{}},	{\Mk{6,2,3,4,5,\bred{1}}{2}}
        },
        {{\Mk{4,5,6,1,2,\bred{3}}{}}
        },
        {{\Mk{5,6,1,2,3,\bred{4}}{}}
        }
        }
	\draw[] (0,0) to (0,6);
}
\qquad\hbox{with}\qquad
\begin{array}{l}
x^T = x_1x_2x_3x_4x_5x_1 \\
\wt(T) = q \frac{(1-t)}{(1-qt^3)}\frac{(1-t)}{(1-qt^4)}\frac{(1-t)}{(1-qt^4)} \\
\mathrm{cov}^z_>(T) = 3+\frac12(9-12-3)=0.
\end{array}
$$
$$
z^{(4)} =
\ \TikZ{[xscale=1.55,yscale=-.6]
     \smTableau{
        {{\Mk{2,3,4,5,6,\bred{1}}{}}, 	{\Mk{6,1,2,3,4,\bred{5}}{1,2}}
        },
        {{\Mk{3,4,5,6,1,\bred{2}}{}},	{\Mk{6,1,3,4,5,\bred{2}}{}}
        },
        {{\Mk{4,5,6,1,2,\bred{3}}{}}
        },
        {{\Mk{5,6,1,2,3,\bred{4}}{}}
        }
        }
	\draw[] (0,0) to (0,6);
}
\qquad\hbox{with}\qquad
\begin{array}{l}
x^T = x_1x_2x_3x_4x_5x_2 \\
\wt(T) = q \frac{(1-t)}{(1-qt^3)}\frac{(1-t)}{(1-qt^4)} \\
\mathrm{cov}^z_>(T) = 3+\frac12(8-12-2)=0.
\end{array}
$$
$$
z^{(5)} =
\ \TikZ{[xscale=1.55,yscale=-.6]
     \smTableau{
        {{\Mk{2,3,4,5,6,\bred{1}}{}}, 	{\Mk{1,6,2,3,4,\bred{5}}{1}}
        },
        {{\Mk{3,4,5,6,1,\bred{2}}{}},	{\Mk{6,2,3,4,5,\bred{1}}{1,2}}
        },
        {{\Mk{4,5,6,1,2,\bred{3}}{}}
        },
        {{\Mk{5,6,1,2,3,\bred{4}}{}}
        }
        }
	\draw[] (0,0) to (0,6);
}
\qquad\hbox{with}\qquad
\begin{array}{l}
x^T = x_1x_2x_3x_4x_5x_1 \\
\wt(T) = q^2t^5 \frac{(1-t)}{(1-qt^4)}\frac{(1-t)}{(1-qt^4)}\frac{(1-t)}{(1-qt^5)} \\
\mathrm{cov}^z_>(T) = 4+4+\frac12(9-12-3)=5.
\end{array}
$$
$$
z^{(6)} =
\ \TikZ{[xscale=1.55,yscale=-.6]
     \smTableau{
        {{\Mk{2,3,4,5,6,\bred{1}}{}}, 	{\Mk{1,6,2,3,4,\bred{5}}{1}}
        },
        {{\Mk{3,4,5,6,1,\bred{2}}{}},	{\Mk{2,6,3,4,5,\bred{1}}{1}}
        },
        {{\Mk{4,5,6,1,2,\bred{3}}{}}
        },
        {{\Mk{5,6,1,2,3,\bred{4}}{}}
        }
        }
	\draw[] (0,0) to (0,6);
}
\qquad\hbox{with}\qquad
\begin{array}{l}
x^T = x_1x_2x_3x_4x_5x_1 \\
\wt(T) = qt \frac{(1-t)}{(1-qt^4)}\frac{(1-t)}{(1-qt^5)} \\
\mathrm{cov}^z_>(T) = 4+\frac12(8-12-2)=1.
\end{array}
$$
$$
z^{(7)} =
\ \TikZ{[xscale=1.55,yscale=-.6]
     \smTableau{
        {{\Mk{2,3,4,5,6,\bred{1}}{}}, 	{\Mk{1,6,2,3,4,\bred{5}}{1}}
        },
        {{\Mk{3,4,5,6,1,\bred{2}}{}},	{\Mk{1,2,3,4,5,\bred{6}}{2}}
        },
        {{\Mk{4,5,6,1,2,\bred{3}}{}}
        },
        {{\Mk{5,6,1,2,3,\bred{4}}{}}
        }
        }
	\draw[] (0,0) to (0,6);
}
\qquad\hbox{with}\qquad
\begin{array}{l}
x^T = x_1x_2x_3x_4x_5x_6 \\
\wt(T) = q^2 t \frac{(1-t)}{(1-qt^4)}\frac{(1-t)}{(1-qt^4)} \\
\mathrm{cov}^z_>(T) = 4+4+\frac12(0-12-2)=1.
\end{array}
$$
$$
z^{(8)} =
\ \TikZ{[xscale=1.55,yscale=-.6]
     \smTableau{
        {{\Mk{2,3,4,5,6,\bred{1}}{}}, 	{\Mk{1,6,2,3,4,\bred{5}}{1}}
        },
        {{\Mk{3,4,5,6,1,\bred{2}}{}},	{\Mk{5,1,3,4,6,\bred{2}}{}}
        },
        {{\Mk{4,5,6,1,2,\bred{3}}{}}
        },
        {{\Mk{5,6,1,2,3,\bred{4}}{}}
        }
        }
	\draw[] (0,0) to (0,6);
}
\qquad\hbox{with}\qquad
\begin{array}{l}
x^T = x_1x_2x_3x_4x_5x_2 \\
\wt(T) = q t \frac{(1-t)}{(1-qt^4)} \\
\mathrm{cov}^z_>(T) = 4+\frac12(7-12-1)=1.
\end{array}
$$
$$
z^{(9)} =
\ \TikZ{[xscale=1.55,yscale=-.6]
     \smTableau{
        {{\Mk{2,3,4,5,6,\bred{1}}{}}, 	{\Mk{5,1,2,3,4,\bred{6}}{2}}
        },
        {{\Mk{3,4,5,6,1,\bred{2}}{}},	{\Mk{1,2,3,4,6,\bred{5}}{1,2}}
        },
        {{\Mk{4,5,6,1,2,\bred{3}}{}}
        },
        {{\Mk{5,6,1,2,3,\bred{4}}{}}
        }
        }
	\draw[] (0,0) to (0,6);
}
\qquad\hbox{with}\qquad
\begin{array}{l}
x^T = x_1x_2x_3x_4x_6x_5 \\
\wt(T) = q^2t \frac{(1-t)}{(1-qt^3)}\frac{(1-t)}{(1-qt^4)}\frac{(1-t)}{(1-qt^5)} \\
\mathrm{cov}^z_>(T) = 3+5+\frac12(1-12-3)=1.
\end{array}
$$
$$
z^{(10)} =
\ \TikZ{[xscale=1.55,yscale=-.6]
     \smTableau{
        {{\Mk{2,3,4,5,6,\bred{1}}{}}, 	{\Mk{5,1,2,3,4,\bred{6}}{2}}
        },
        {{\Mk{3,4,5,6,1,\bred{2}}{}},	{\Mk{2,1,3,4,6,\bred{5}}{1}}
        },
        {{\Mk{4,5,6,1,2,\bred{3}}{}}
        },
        {{\Mk{5,6,1,2,3,\bred{4}}{}}
        }
        }
	\draw[] (0,0) to (0,6);
}
\qquad\hbox{with}\qquad
\begin{array}{l}
x^T = x_1x_2x_3x_4x_6x_5 \\
\wt(T) = q^2t^2 \frac{(1-t)}{(1-qt^3)}\frac{(1-t)}{(1-qt^5)} \\
\mathrm{cov}^z_>(T) = 3+5+\frac12(2-12-2)=2.
\end{array}
$$
$$
z^{(11)} =
\ \TikZ{[xscale=1.55,yscale=-.6]
     \smTableau{
        {{\Mk{2,3,4,5,6,\bred{1}}{}}, 	{\Mk{5,1,2,3,4,\bred{6}}{2}}
        },
        {{\Mk{3,4,5,6,1,\bred{2}}{}},	{\Mk{5,2,3,4,6,\bred{1}}{2}}
        },
        {{\Mk{4,5,6,1,2,\bred{3}}{}}
        },
        {{\Mk{5,6,1,2,3,\bred{4}}{}}
        }
        }
	\draw[] (0,0) to (0,6);
}
\qquad\hbox{with}\qquad
\begin{array}{l}
x^T = x_1x_2x_3x_4x_6x_1 \\
\wt(T) = q \frac{(1-t)}{(1-qt^3)}\frac{(1-t)}{(1-qt^4)} \\
\mathrm{cov}^z_>(T) = 3+\frac12(8-12-2)=0.
\end{array}
$$
$$
z^{(12)} =
\ \TikZ{[xscale=1.55,yscale=-.6]
     \smTableau{
        {{\Mk{2,3,4,5,6,\bred{1}}{}}, 	{\Mk{5,1,2,3,4,\bred{6}}{2}}
        },
        {{\Mk{3,4,5,6,1,\bred{2}}{}},	{\Mk{5,1,3,4,6,\bred{2}}{}}
        },
        {{\Mk{4,5,6,1,2,\bred{3}}{}}
        },
        {{\Mk{5,6,1,2,3,\bred{4}}{}}
        }
        }
	\draw[] (0,0) to (0,6);
}
\qquad\hbox{with}\qquad
\begin{array}{l}
x^T = x_1x_2x_3x_4x_6x_2 \\
\wt(T) = q \frac{(1-t)}{(1-qt^3)} \\
\mathrm{cov}^z_>(T) = 3+\frac12(7-12-1)=0.
\end{array}
$$
$$
z^{(13)} =
\ \TikZ{[xscale=1.55,yscale=-.6]
     \smTableau{
        {{\Mk{2,3,4,5,6,\bred{1}}{}}, 	{\Mk{5,6,2,3,4,\bred{1}}{}}
        },
        {{\Mk{3,4,5,6,1,\bred{2}}{}},	{\Mk{6,2,3,4,1,\bred{5}}{1,2}}
        },
        {{\Mk{4,5,6,1,2,\bred{3}}{}}
        },
        {{\Mk{5,6,1,2,3,\bred{4}}{}}
        }
        }
	\draw[] (0,0) to (0,6);
}
\qquad\hbox{with}\qquad
\begin{array}{l}
x^T = x_1x_2x_3x_4x_1x_5 \\
\wt(T) = qt \frac{(1-t)}{(1-qt^4)}\frac{(1-t)}{(1-qt^5)} \\
\mathrm{cov}^z_>(T) = 4+\frac12(8-12-2)=1.
\end{array}
$$
$$
z^{(14)} =
\ \TikZ{[xscale=1.55,yscale=-.6]
     \smTableau{
        {{\Mk{2,3,4,5,6,\bred{1}}{}}, 	{\Mk{5,6,2,3,4,\bred{1}}{}}
        },
        {{\Mk{3,4,5,6,1,\bred{2}}{}},	{\Mk{2,6,3,4,1,\bred{5}}{1}}
        },
        {{\Mk{4,5,6,1,2,\bred{3}}{}}
        },
        {{\Mk{5,6,1,2,3,\bred{4}}{}}
        }
        }
	\draw[] (0,0) to (0,6);
}
\qquad\hbox{with}\qquad
\begin{array}{l}
x^T = x_1x_2x_3x_4x_6x_5 \\
\wt(T) = qt^2 \frac{(1-t)}{(1-qt^5)} \\
\mathrm{cov}^z_>(T) = 5+\frac12(7-12-1)=2.
\end{array}
$$
$$
z^{(15)} =
\ \TikZ{[xscale=1.55,yscale=-.6]
     \smTableau{
        {{\Mk{2,3,4,5,6,\bred{1}}{}}, 	{\Mk{5,6,2,3,4,\bred{1}}{}}
        },
        {{\Mk{3,4,5,6,1,\bred{2}}{}},	{\Mk{5,2,3,4,1,\bred{6}}{2}}
        },
        {{\Mk{4,5,6,1,2,\bred{3}}{}}
        },
        {{\Mk{5,6,1,2,3,\bred{4}}{}}
        }
        }
	\draw[] (0,0) to (0,6);
}
\qquad\hbox{with}\qquad
\begin{array}{l}
x^T = x_1x_2x_3x_4x_6x_1 \\
\wt(T) = qt \frac{(1-t)}{(1-qt^4)} \\
\mathrm{cov}^z_>(T) = 4+\frac12(7-12-1)=1.
\end{array}
$$
$$
z^{(16)} =
\ \TikZ{[xscale=1.55,yscale=-.6]
     \smTableau{
        {{\Mk{2,3,4,5,6,\bred{1}}{}}, 	{\Mk{5,6,2,3,4,\bred{1}}{}}
        },
        {{\Mk{3,4,5,6,1,\bred{2}}{}},	{\Mk{5,6,3,4,1,\bred{2}}{}}
        },
        {{\Mk{4,5,6,1,2,\bred{3}}{}}
        },
        {{\Mk{5,6,1,2,3,\bred{4}}{}}
        }
        }
	\draw[] (0,0) to (0,6);
}
\qquad\hbox{with}\qquad
\begin{array}{l}
x^T = x_1x_2x_3x_4x_1x_2 \\
\wt(T) = 1 \\
\mathrm{cov}^z_>(T) = 0+\frac12(12-12-0)=0.
\end{array}
$$
In summary,
\begin{align*}
(x_1x_2x_3x_4)^{-1}&E_{(2,2,1,1,0,0)}\\
=q^3&\frac{(1-t)}{(1-qt^3)}\frac{(1-t)}{(1-qt^4)}\frac{(1-t)}{(1-qt^4)}\frac{(1-t)}{(1-qt^5)} x_5x_6
+q^2 t \frac{(1-t)}{(1-qt^3)}\frac{(1-t)}{(1-qt^4)}\frac{(1-t)}{(1-qt^5)} x_5x_6
\\
&+q \frac{(1-t)}{(1-qt^3)}\frac{(1-t)}{(1-qt^4)}\frac{(1-t)}{(1-qt^4)} x_5x_1
+q \frac{(1-t)}{(1-qt^3)}\frac{(1-t)}{(1-qt^4)} x_5x_2 
\\
&+q^2t^5 \frac{(1-t)}{(1-qt^4)}\frac{(1-t)}{(1-qt^4)}\frac{(1-t)}{(1-qt^5)} x_5x_1
+qt \frac{(1-t)}{(1-qt^4)}\frac{(1-t)}{(1-qt^5)} x_5x_1
\\
&+q^2 t \frac{(1-t)}{(1-qt^4)}\frac{(1-t)}{(1-qt^4)} x_5x_6
+q t \frac{(1-t)}{(1-qt^4)}  x_5x_2
\\
&+q^2t \frac{(1-t)}{(1-qt^3)}\frac{(1-t)}{(1-qt^4)}\frac{(1-t)}{(1-qt^5)} x_6x_5
+q^2t^2 \frac{(1-t)}{(1-qt^3)}\frac{(1-t)}{(1-qt^5)} x_6x_5
\\
&+q \frac{(1-t)}{(1-qt^3)}\frac{(1-t)}{(1-qt^4)} x_6x_1
+q \frac{(1-t)}{(1-qt^3)}  x_6x_2
\\
&+qt \frac{(1-t)}{(1-qt^4)}\frac{(1-t)}{(1-qt^5)} x_1x_5
+qt^2 \frac{(1-t)}{(1-qt^5)} x_6x_5
+qt \frac{(1-t)}{(1-qt^4)} x_6x_1
+x_1x_2.
\end{align*}

\bigskip\noindent
\textbf{Acknowledgements.}  We thank A.\ Garbali and Tianyi Yu for very helpful proofreading and corrections.

\end{document}